\documentclass[12pt]{article}

\usepackage{epsfig}
\usepackage{fullpage}
\usepackage{amsfonts}
\usepackage{amsmath}
\usepackage{amsthm}
\newtheorem{theorem}{Theorem}
\newtheorem{lemma}{Lemma}
\usepackage{url}

\newcommand{\HH}{{\mathbb H}}

\newcommand{\OO}{{\mathbb O}}
\newcommand{\RR}{{\mathbb R}}
\renewcommand{\Re}{{\rm Re}}
\renewcommand{\Im}{{\rm Im}}
\renewcommand{\bar}{\overline}

\begin{document}

\title{\boldmath\textbf{Finding $3\times3$ Hermitian Matrices over the
Octonions with Imaginary Eigenvalues}}

\author{
	Henry Gillow-Wiles\\
	Department of Science and Mathematics Education\\
	Oregon State University\\
	Corvallis, OR  97331\\
	\texttt{gillowwh{\rm @}onid.orst.edu}\\
	\and
	Tevian Dray\\
	Department of Mathematics\\
	Oregon State University\\
	Corvallis, OR  97331\\
	\texttt{tevian{\rm @}math.oregonstate.edu}\\
	}

\maketitle

\begin{abstract}
We show that any $3$-component octonionic vector which is purely imaginary,
but not quaternionic, is an eigenvector of a 6-parameter family of Hermitian
octonionic matrices, with imaginary eigenvalue equal to the associator of its
elements.
\end{abstract}

\section{Introduction}

The eigenvalue problem for $3\times3$ Hermitian octonionic matrices,
henceforth referred to as \textit{Jordan matrices}, contains some surprises.
Notable among these is that, whereas each Jordan matrix satisfies its
characteristic equation, its real eigenvalues do not.  As shown
in~\cite{Eigen,Bases}, each Jordan matrix admits six real eigenvalues, rather
than three.  However, the eigenvalues divide naturally into two families of
three, and the corresponding families of eigenvectors do have the expected
properties, such as orthonormality, provided that these properties are
properly formulated.

Due to the nonassociativity of the octonions, most Jordan matrices also appear
to admit eigenvalues which are not real; several examples were discussed
in~\cite{NonReal,Properties}.  However, to our knowledge there is no known
algorithm for finding the non-real eigenvalues of such matrices, nor is it
clear how many there are.

In this paper, which is an extension of~\cite{HenryThesis}, we take a
different approach.  Rather than attempt to find the non-real eigenvalues and
corresponding eigenvectors of a given Jordan matrix, we instead find the
Jordan matrices which admit a given eigenvector/eigenvalue pair.
Specifically, for any vector $v\in\OO^3$ which is not quaternionic, we use the
associator of the elements of $v$ as the eigenvalue, and find all Jordan
matrices for which $v$ is an eigenvector with that eigenvalue, which is
nonzero by assumption.  We show below that a necessary condition for such
matrices to exist is that $\Re(v)=0$, and that if this condition is satisfied
there is a 6-parameter family of such matrices.

We begin in Section~\ref{octonions} by reviewing the octonions and their
properties, and then briefly summarize some known examples~\cite{Properties}
of $3\times3$ Hermitian octonionic matrices with imaginary eigenvalues in
Section~\ref{examples}.  In Section~\ref{results}, we present our new results,
which we then summarize in Section~\ref{conclusions}, where we also propose
some further conjectures.

\section{Octonions}
\label{octonions}

\begin{figure}
\begin{center}
\epsfxsize=3.5in
\epsffile[68 168 543 614]{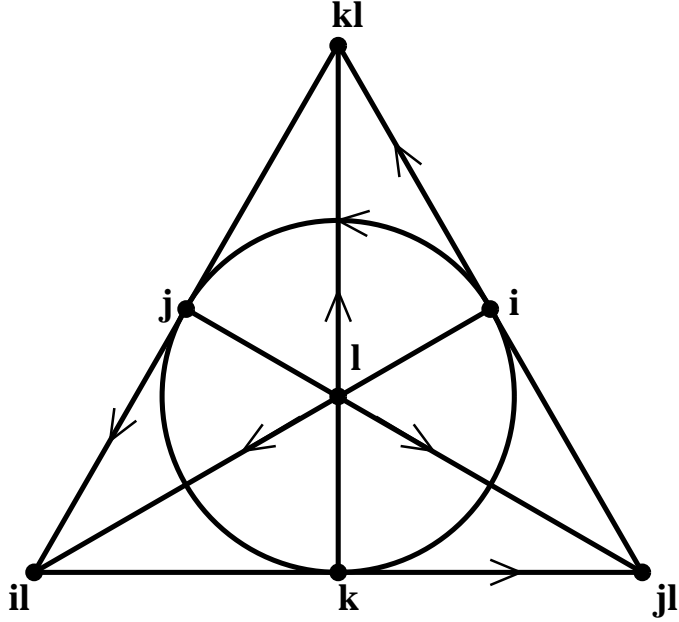}
\end{center}
\caption{The representation of the octonionic multiplication table using the
7-point projective plane.  Each of the 7 oriented lines gives a
quaternionic triple.}
\label{omult}
\end{figure}

We use the standard basis $\{1,i,j,k\}$ for the quaternions $\HH$, and we
construct the octonions $\OO$ via the Cayley-Dickson process as
$\HH\oplus\HH\ell$.  The resulting multiplication table is neatly summarized
by the oriented Fano geometry shown in Figure~\ref{omult}.  As is well-known,
the octonions are neither commutative nor associative.

Writing the components of an octonion $w$ as
\begin{equation}
w = w_1 + w_2 \,i + w_3 \,j + w_4 \,k
	+ w_5 \,k\ell + w_6 \,j\ell + w_7 \,i\ell + w_8 \,\ell
\end{equation}
we have
\begin{eqnarray}
\Re(w) &=& w_1 \\
\Im(w) &=& w - \Re(w)
\end{eqnarray}
and
\begin{eqnarray}
\bar{w} &=& 2\,\Re(w) - w \\
|w|^2 &=& w\,\bar{w}
\end{eqnarray}
Any three octonions $x,y,z\in\OO$ can be assumed without loss of generality to
take the form
\begin{eqnarray}
x &=& x_1 + x_2 \,i \\
y &=& y_1 + y_2 \,i + y_3 \,j \\
z &=& z_1 + z_2 \,i + z_3 \,j + z_4 \,k + z_8 \,\ell
\end{eqnarray}
by suitable choice of basis; we refer to $x,y,z$ as \textit{generic
octonions}.  Choosing
\begin{equation}
v =
\begin{pmatrix}
x\cr y\cr z\cr
\end{pmatrix}
\in\OO^3
\end{equation}
the associator of $v$ is by definition the associator of its elements, that is
\begin{equation}
[v] = [x,y,z] = (xy)z - x(yz)
    = 2 x_2 y_3 z_8 \,k\ell
\end{equation}
with the last equality holding for generic octonions.

We seek solutions of the (right) eigenvalue problem
\begin{equation}
A v = v \lambda
\label{eigen}
\end{equation}
where $\lambda\not\in\RR$; in what follows we will consider only the case
$\lambda = [v] \ne 0$, noting that $[v]$ is pure imaginary.  Any Jordan matrix
can be written in the form
\begin{equation}
A =
\begin{pmatrix}
p& a& \bar{c}\cr
\bar{a}& m& b\cr
c& \bar{b}& n\cr
\end{pmatrix}
\end{equation}
with $p,m,n\in\RR$ and $a,b,c\in\OO$.  Then~\eqref{eigen} takes the form
\begin{eqnarray}
p x + a y + \bar{c} z &=& x \lambda \label{I}\\
\bar{a} x + m y + b z &=& y \lambda \label{II}\\
c x + \bar{b} y + n z &=& z \lambda \label{III}
\end{eqnarray}

As shown below,~\eqref{eigen} admits solutions only if $\Re(v)=0$, in which
case there is a \hbox{6-parameter} family of Jordan matrices $A$ which satisfy
it.

\section{Examples}
\label{examples}

Little is known about solutions of~\eqref{eigen} with non-real eigenvalues.
Although the problem of finding real eigenvalues for $3\times3$ Hermitian
octonionic matrices has been completely solved~\cite{Eigen}, no such solution
exists for finding imaginary eigenvalues~\cite{NonReal,Properties}.  In fact,
we know of only a handful of explicit examples of families of $3\times3$
Hermitian octonionic matrices admitting imaginary eigenvalues, such as those
given in~\cite{Properties}, which are reproduced below.  Note in each case
that $p$ can always be chosen so that the eigenvalue has no real part.
Furthermore, in the first example, the eigenvectors have no real part, the
associator of each eigenvector is a (possibly zero) multiple of $k\ell$, the
imaginary direction of the eigenvalue, and this last property also holds for
the elements of $A_1$.  However, the first two properties fail to hold in the
second example, while the last property fails in the third.

\paragraph{Example 1}

The matrix
\begin{equation}
A_1 =
\begin{pmatrix}
~~p & ~~iq & ~kqs\\
-iq &~~ p & ~jq\\
-kqs & -jq & ~p
\end{pmatrix}
\label{ex1}
\end{equation}
with $p,q\in\RR$ and
\begin{equation}
s = \cos \theta +k\ell \sin \theta
\end{equation}
has, among others, the eigenvalues and eigenvectors,
\begin{subequations}
\begin{align}
\lambda_{u} &= p\pm q\bar{s}: &
  u _\pm &= \begin{pmatrix}i\\0\\j\end{pmatrix}S_ \pm\\
\lambda_{v} &= p\pm q\bar{s}: &
  v_ \pm &= \begin{pmatrix}j\\ 2ks\\i\end{pmatrix}S_ \pm\\
\lambda_{w} &= p\mp 2q\bar{s}: &
  w_ \pm &= \begin{pmatrix}j\\-ks\\i\end{pmatrix}S_ \pm
\end{align}
\end{subequations}
where
\begin{equation}
S_ \pm = \left\{ \begin{matrix}-k\,\ell\\1\end{matrix} \right.
\end{equation}

\paragraph{Example 2}

The matrix
\begin{equation}
A_2 =
\begin{pmatrix}
~~p & ~~qi & \frac{q}{6}(\sqrt{5}k+2 \ell)\\
\noalign{\smallskip}
-qi &~~ p & ~\frac{q}{2}j\\
\noalign{\smallskip}
- \frac{q}{6}(\sqrt{5}k+2 \ell) & -\frac{q}{2}j & ~p
\end{pmatrix}
\label{ex2}
\end{equation}
has, among others, the eigenvectors and
eigenvalues,
\begin{subequations}
\begin{align}
\lambda_{u_1} &= (p+\frac{\sqrt{5}}{2}q)-\frac{q}{2}k\ell: &
  u_1 &= \begin{pmatrix}
	3k\\ \sqrt{5}j-2\,i\ell\\1+ \sqrt{5}\,k\ell
	\end{pmatrix}\\
\lambda_{u_2} &= (p+\frac{\sqrt{5}}{2}q)+\frac{q}{2}k\ell: &
  u_2 &= \begin{pmatrix}\sqrt{5}k+2 \ell\\ 3j\\ \sqrt{5}-k\ell\end{pmatrix}\\
\lambda_{v_1} &= (p-\frac{\sqrt{5}}{3}q)+\frac{2q}{3}k\ell:&
  v_1 &= \begin{pmatrix}\sqrt{5}j-2\,i\ell\\ 3k\\ 0\end{pmatrix}\\
\lambda_{v_2} &= (p-\frac{\sqrt{5}}{3}q)-\frac{2q}{3}k\ell: &
  v_2 &= \begin{pmatrix}3j \\ \sqrt{5}k+2\ell\\0\end{pmatrix}\\
\lambda_{w_1} &= (p-\frac{\sqrt{5}}{6}q)-\frac{q}{6}k\ell:&
  w_1 &= \begin{pmatrix}3k\\\sqrt{5}j-2\,i\ell\\ -7-\sqrt{5}k\ell\end{pmatrix}\\
\lambda_{w_2} &= (p-\frac{\sqrt{5}}{6}q)+\frac{q}{6}k\ell:&
  w_2 &= \begin{pmatrix}\sqrt{5}k+2 \ell\\3j\\-3\sqrt{5}-3\,k\ell\end{pmatrix}
\end{align}
\end{subequations}

\paragraph{Example 3}

The matrix
\begin{equation}
A_3 =
\begin{pmatrix}
~~p & ~~qi & -q(j-i\ell-j\ell)\\
\noalign{\smallskip}
-qi &~~ p & q(1+k+\ell)\\
\noalign{\smallskip}
q(j-i\ell-j\ell) & q(1-k-\ell) & ~p
\end{pmatrix}
\label{ex3}
\end{equation}
admits the eigenvector
\begin{equation}
v = \begin{pmatrix}j\\\ell\\0\end{pmatrix}
\end{equation}
with eigenvalue
\begin{equation}
\lambda_v=p-q k\ell
\end{equation}

\goodbreak
\newpage

\section{Results}
\label{results}

As already noted, for each example in the previous section, $p$ can be chosen
so that a given eigenvalue is purely imaginary.  More generally, the real part
of an eigenvalue can be changed by adding a suitable multiple of the identity
matrix to the original matrix.  More formally, we have the following result:

\begin{lemma}
The Hermitian matrices having a given eigenvector can be divided into families
which differ only by (real) multiples of the identity matrix.  The imaginary
part of the corresponding eigenvalue is the same for each member of such a
family, and each family contains a unique member such that the real part of
the corresponding eigenvalue vanishes.
\label{nonreal}
\end{lemma}

\begin{proof}
If $v$ is an eigenvector of $A$ with eigenvalue $\lambda$, then $v$ is also an
eigenvector of $A+pI$ for any $p\in\RR$, with eigenvalue $\lambda+p$.  In
particular, $v$ is an eigenvector of $A-\Re(\lambda)I$, with eigenvalue
$\Im(\lambda)$.
\end{proof}

Note that this technique can in general be used to eliminate the real part of
only one eigenvalue at a time.  Nonetheless, any eigenvector with a non-real
eigenvalue is also an eigenvector of a closely related matrix with a purely
imaginary eigenvalue.

This suggests the following strategy for trying to find eigenvectors with
non-real eigenvalues: Rather than trying to find eigenvector/eigenvalue pairs
$v$ and $\lambda$ satisfying~\eqref{eigen} for given $A$, with $\lambda$
non-real, we will instead seek to categorize the matrices which admit such
eigenvalues.  Lemma~\ref{nonreal} now tells us that we can assume
$\Re(\lambda$)=0 without loss of generality, at least so long as we consider
only a single eigenvector.  We will therefore attempt to find the matrices $A$
which admit a given vector $v$ as an eigenvector, with given eigenvalue
$\lambda$ satisfying $\Re(\lambda$)=0.  Motivated by the first example, we
will further assume that $\lambda$ is a real multiple of $[v]$, and we will
consider only the case where $[v]\ne0$.  Finally, by rescaling $A$, the
constant of proportionality can be assumed to be $1$.  Thus, we assume that
\begin{equation}
\lambda = [v] \ne 0
\end{equation}

Since $[v]\ne0$ by assumption, none of $x,y,z$ can be zero.  In particular,
$x\ne0$, and it is straightforward to solve~\eqref{II} for $\bar{a}$
and~\eqref{III} for $c$, yielding
\begin{eqnarray}
\bar{a} &=& \Bigl({y(\lambda-m) - bz}\Bigr) \frac{\bar{x}}{~|x|^2} \\
c &=& \Bigl({z(\lambda-n) - \bar{b}y}\Bigr) \frac{\bar{x}}{~|x|^2}
\end{eqnarray}
Inserting these expressions into~\eqref{I} reduces~\eqref{eigen} to the form
\begin{equation}
\Bigl(x \left(\bar{\lambda}\bar{y}-\bar{z}\bar{b}\right) \Bigr) y
	+ \Bigl(x \left(\bar{\lambda}\bar{z}-\bar{y}b\right) \Bigr) z
	- x \lambda |x|^2
  = x \left( m|y|^2 + n|z|^2 - p|x|^2 \right)
\label{Master}
\end{equation}

\goodbreak

\begin{lemma}
\textit{If $[v]\ne0$, then $b\perp\lambda$, that is, $b_5=0$.}
\end{lemma}

\begin{proof}
Multiply both sides of~\eqref{Master} on the left by $\bar{x}$.  The LHS of
the resulting expression must be real, since the RHS is, but direct
computation shows that the coefficient of $i$ on the left is
$2 |x|^2 y_3 z_8 b_5$.  Since each factor except for $b_5$ is nonzero by
assumption, $b_5$ must be zero.
\end{proof}

\goodbreak

\begin{theorem}
If $[v]\ne0$, then there are no solutions to~\eqref{eigen} unless
$\Re(v)=0$.
\end{theorem}

\begin{proof}
Direct computation, as follows.  Inserting $b_5=0$ into~\eqref{Master}, the
$j$-component yields $4 x_2^2 y_3 z_8^2 z_1 = 0$.  Since each factor except
for $z_1$ is nonzero by assumption, $\Re(z)=0$.  In a separate computation,
the $k$-component of~\eqref{Master} can be solved for $b_8$, yielding
\[
b_8 = \frac{2 y_3 z_2 z_8 {x_2}^2 + b_6 y_1 x_2 - b_6 x_1 y_2 + b_7 x_1 y_3}
	{x_2 y_3}
\]
Inserting the
result into the $\ell$-component of~\eqref{Master}, along with $b_5=0=z_1$,
results in $-4 x_2^2 y_3^2 z_8 y_1 = 0$, which forces $\Re(y)=0$.  Finally,
the $i\ell$- and $j\ell$-components of~\eqref{Master} can be solved for $b_2$
and $b_3$, yielding
\begin{eqnarray*}
b_2 &=& \frac{1}{{x_2}^2 {y_3}^2 z_8} \big( 
	{y_3}^2 z_8 {x_2}^5 - {y_3}^2 {z_8}^3 {x_2}^3 - {y_3}^4 z_8 {x_2}^3
	+ {x_1}^2 {y_3}^2 z_8 {x_2}^3 + {y_2}^2 {y_3}^2 z_8 {x_2}^3 \\
&&	+ {y_3}^2 {z_2}^2 z_8 {x_2}^3 - {y_3}^2 {z_3}^2 z_8 {x_2}^3
	- {y_3}^2 {z_4}^2 z_8 {x_2}^3 + b_6 {y_3}^2 z_4 {x_2}^2
	- 2 {x1} y_2 y_3 z_2 z_4 z_8 {x_2}^2 \\
&&	+ b_7 x_1 {y_3}^2 z_2 x_2 - b_7 x_1 y_2 y_3 z_3 x_2
	+ b_4 x_1 y_2 y_3 z_8 x_2 + b_6 {x_1}^2 {y_2}^2 z_4
	- b_7 {x_1}^2 y_2 y_3 z_4 \big) \\
b_3 &=& \frac{1} {{x_2}^2 y_3 z_8} \big(
	2 y_2 {y_3}^2 z_8 {x_2}^3 + 2 y_3 z_2 z_3 z_8 {x_2}^3
	- b_7 y_3 z_4 {x_2}^2 - 2 x_1 y_3 z_2 z_4 z_8 {x_2}^2 \\
&&	+ b_6 x_1 y_3 z_2 x_2 - b_6 x_1 y_2 z_3 x_2
	+ b_4 x_1 y_3 z_8 x_2 + b_6 {x_1}^2 y_2 z_4
	- b_7 {x_1}^2 y_3 z_4 \big)
\end{eqnarray*}
and the result inserted into the $k$-component (along with $b_5=0=z_1=y_1$ and
the above expression for $b_8$), resulting in $-2 x_1 |x|^2 \lambda = 0$,
which forces $\Re(x)=0$.  Thus, $\Re(v)=0$.
\end{proof}

\begin{theorem}
\textit{If $[v]\ne0$ and $\Re(v)=0$, then there is a 6-parameter family of
solutions to~\eqref{eigen}}.
\label{result}
\end{theorem}

\begin{proof}
Inserting the above expressions for $b_2$, $b_3$ and $b_8$, as well as the
condition $\Re(v)=0$, into~\eqref{Master} results in a single nonzero
component, which can be solved for $b_6$, yielding
\begin{eqnarray*}
b_6 &=& \big(
	- 2 y_3 z_4 z_8 {x_2}^3 - p z_8 {x_2}^2
	+ 2 y_3 z_4 {z_8}^3 x_2 + 2 y_3 {z_4}^3 z_8 x_2
	+ 2 {y_3}^3 z_4 z_8 x_2 - 2 y_3 {z_2}^2 z_4 z_8 x_2 \\
&&	+ 2 y_3 {z_3}^2 z_4 z_8 x_2 + 2 {y_2}^2 y_3 z_4 z_8 x_2
	+ 4 y_2 z_2 z_3 z_4 z_8 x_2 + n {z_8}^3 - 2 b_7 y_2 {z_4}^2
	- 2 b_7 y_2 {z_8}^2 \\
&&	+ m {y_2}^2 z_8 + m {y_3}^2 z_8
	+ n {z_2}^2 z_8 + n {z_3}^2 z_8 + n {z_4}^2 z_8 \\
&&	+ 2 b_1 y_2 z_2 z_8 + 2 b_4 y_3 z_2 z_8 - 2 b_4 y_2 z_3 z_8
	+ 2 b_1 y_3 z_3 z_8\big)
	\big/\big(2 y_3 \big({z_4}^2 + {z_8}^2\big)\big)
\end{eqnarray*}
As with the equations solved above for $b_2$, $b_3$ and $b_8$, the relevant
coefficients are nonzero under the stated assumptions, so that the given
solutions always exist.  We have thus constructed $A$ explicitly, with $b_1$,
$b_4$, $b_7$, $p$, $m$, and $n$ as free parameters.
\end{proof}

\section{Conclusion}
\label{conclusions}

We have created a method for finding a Hermitian matrix $A\in\OO^{3\times3}$
which has an eigenvalue relationship with an imaginary vector $v\in\OO^3$,
with the associator of $v$, assumed to be nonzero, playing the role of the
eigenvalue $\lambda$.  For our method to be successful, rather than fixing the
matrix, we must begin by fixing the vector $v$, thus fixing $\lambda$ as well.

The question we must ask is if the resulting eigenvalue/eigenvector system
from our construction method represents a variation of one of the existing
three family examples presented by Dray, Janesky and
Manogue~\cite{Properties}.

The last two eigenvectors in Example~1 satisfy the conditions of
Theorem~\ref{result} (vanishing real part and non-vanishing associator), and
$p$ and $q$ can be chosen so that the eigenvalue is precisely the vector
associator, as required by our hypotheses.  It is straightforward to verify
that the corresponding matrix $A_1$ is indeed contained in the 6-parameter
family constructed in Theorem~\ref{result}.  Lemma~\ref{nonreal} can now be
used, along with some obvious renormalization, to construct the matrices $A_1$
for any values of $p$ and $q$.  In this sense, Example 1 is contained within
our solution method, although our method generates many more solutions --- but
can not yet handle the first eigenvector shown, whose associator vanishes.

Each of the eigenvectors in Examples~2 and~3, however, either has a non-zero
real part or a vanishing associator, so that our results do not apply to these
cases.  Note that in Example~2, although the imaginary part of the eigenvalues
is indeed in the direction of $[A_2]$ (the associator of the elements of
$A_2$), namely $k\ell$, none of the given eigenvectors has an associator in
this direction.  Furthermore, in Example~3, the imaginary part of the
eigenvalue no longer points in the direction of the matrix associator $[A_3]$.
These examples therefore make clear that the results in this paper represent
only the tip of the iceberg; our assumptions are too restrictive.

By understanding more about the problems encountered in trying to find a
characteristic eigenvalue equation for $3\times3$ Hermitian matrices over the
octonions, we hope our work will aid in the discovery of a method for finding
the imaginary eigenvalues (if any), and their corresponding eigenvectors, for
any given $3\times3$ Hermitian octonionic matrix.

%\bibliographystyle{TDbib}
%\bibliography{octonions,OctoRefs}

\newpage

\begin{thebibliography}{1}

\bibitem{Eigen}
Tevian Dray and Corinne~A. Manogue, {\em {The Octonionic Eigenvalue Problem}},
  Adv. Appl. Clifford Algebras {\bf 8}, 341--364 (1998).

\bibitem{Bases}
Tevian Dray, Corinne~A. Manogue, and Susumu Okubo, {\em {Orthonormal Eigenbases
  over the Octonions}}, Algebras Groups Geom. {\bf 19}, 163--180 (2002).

\bibitem{NonReal}
Tevian Dray, Jason Janesky, and Corinne~A. Manogue, {\em {Octonionic Hermitian
  Matrices with Non-Real Eigenvalues}}, Adv. Appl. Clifford Algebras {\bf 10},
  193--216 (2000).

\bibitem{Properties}
Tevian Dray, Jason Janesky, and Corinne~A. Manogue,
\newblock {\it Some Properties of $3\times3$ Octonionic Hermitian Matrices with
  Non-Real Eigenvalues},
\newblock Technical report, Oregon State University (2000).
\newblock \null\hfill(\url{http://xxx.lanl.gov/abs/math/0010255})

\bibitem{HenryThesis}
Henry Gillow-Wiles,
\newblock {\it Finding $3\times3$ Hermitian Matrices over the Octonions with
  Imaginary Eigenvalues},
\newblock Master's thesis, Oregon State University (2008).

\bibitem{Baez}
John~C. Baez, {\em \it The Octonions}, Bull. Amer. Math. Soc. {\bf 39},
  145--205 (2002).

\end{thebibliography}

\end{document}